\documentclass[11pt]{article}
\usepackage{amsmath}
\usepackage{amsfonts}
\usepackage{graphicx}
\usepackage[sort,compress]{cite}

\begin{document}
\begin{center}
\Large\bf{Estimating the uncertainty in underresolved nonlinear dynamics}
\end{center}

\begin{center}
Alexandre J. Chorin$^{1,2}$, Ole H. Hald$^{1,2}$ 
\vspace{3mm}

$^1$Department of Mathematics\\
University of California, Berkeley\\
Berkeley CA 94720;
\vspace{1mm}

$^2$Lawrence Berkeley National Laboratory.\\
Berkeley CA 94720.
\vspace{1mm}

\end{center}
\vspace{.5mm}
\begin{center}
\emph{Abstract}
\end{center}
The Mori-Zwanzig formalism of statistical mechanics is 
used to estimate the uncertainty caused by underresolution in the solution of 
a nonlinear dynamical system. A general approach is outlined and 
applied to a simple example. The
noise term that describes the uncertainty turns out to be neither
Markovian nor Gaussian. It is argued that this is the general situation.

\section{Introduction}

There are many problems in science 
where it is impractical to find fully resolved numerical solutions of the problems of interest, and one needs to estimate
the resulting uncertainty. For example, there are situations, e.g. in geophysics and in economics, where one wants to estimate the state of a system
on the basis of noisy underresolved equations supplemented by a stream of noisy
data (see e.g. \cite{Douc,Ch3}); this is the ``filtering" or ``data assimilation" problem, which requires that a probability density of the noise
be provided. There are similar problems in turbulence calculations, where 
one wants to model the impact of the motion at one set of scales on the motion at another scale 
(see e.g \cite{Bens,Kuksin}), i.e. one wants to estimate the difference between a full system
and one in which a given set of components is absent. In such settings, one often makes 
the assumption that the missing information can be represented as 
white noise, though it
is well understood that this is justifiable only where the scales are separated, 
i.e., when the components that are
not computed are faster and smaller than the ones that one is solving for (see e.g. \cite{Papa}), and that separation of scales fails to hold in many problems of interest.  
Here 
we present a construction that does not rely on separation of scale, and use it to
illustrate the intrinsic features of the problem. 

It is important to note that
what one needs is an estimate of the uncertainty
in the equations, not in the results. For example, in data assimilation, it is the uncertainty in
the equations that is combined with the uncertainty in the data to produce an 
estimate of the overall uncertainty in the state. Consider for example  
a system of ordinary differential equations of the
form 
\begin{equation}
\label{main}
\frac {d}{dt}\phi=R(\phi),
\end{equation}
with initial data $\phi(0)=x$. Assume that $\phi$ has components $\phi_1,\phi_2,\dots,\phi_K$, only the first $k < K$ of which can be effectively computed.
Define ${\hat \phi}=(\phi_1,\dots,\phi_k)$ and ${\tilde \phi}=(\phi_{k+1},\dots,\phi_K)$. Similarly, partition the right-hand side of the equation as 
$R=(\hat{R},\tilde{R})$ and the data as $x=(\hat{x},\tilde{x})$. The underesolved equations that
we are able to solve have the form
\begin{equation}
\label{psi0}
\frac {d}{dt} \psi =\mathcal{R}(\psi),
\end{equation}
where $\psi$ is an approximation of $\hat{\phi}$, $\mathcal{R}$ is an approximation of $\hat{R}$, and $\psi(0)=\hat{x}$ (when one underresolves, one approximates functions of
many variables by functions of fewer variables, as when one uses $k$ Fourier components to
approximate a differential operator when $K > k$ are needed). 
The estimate of the uncertainty is a function $n(t)$, the ``noise", which,
added to $\mathcal{R}$, would make the solution of equation (\ref{psi0}) approximate the first $k$ components of the solution of equation (\ref{main}); i.e.,
we are hoping that approximately, in a sense to be determined, one would have:
\begin{equation}
\frac {d}{dt}\hat{\phi}=\mathcal{R}(\hat{\phi})+n(t),
\label{psi}
\end{equation}
with initial data $\phi (0)=\hat{x}$. The term $n(t)$ would be the estimate of the uncertainty in the underresolved system (\ref{psi0}). Equation (\ref{psi}) obviously holds if one sets
\begin{equation}
n(t)=\hat{R}(\phi)-\mathcal{R}(\hat{\phi}). 
\label{trick}
\end{equation}
So that, if one can find $\phi$, one can estimate $n(t)$. But if one has $\phi$, there is no need to find $\psi$ or $n$.

The problem becomes more interesting if it is randomized as in \cite{Ford,Zw3} and in the 
Mori-Zwanzig theory (see e.g. \cite{Eva1,Fic,Zw4,Ch13,Darve}). 
Equations (\ref{main}) and equations (\ref{psi0}) need initial data
of different dimension; to compare their solutions one has to decide what to
do with the extra components $\tilde{x}$ of the initial vector $x$.
In practical problems $\tilde{x}$ is unknown, and 
we assume that it is sampled 
from some suitable known probability density, maybe deduced from previous
knowledge, for example, if the problem comes from weather forecasting, ample
records of previous weathers can be put to use. 
The function $\phi(t)$  
becomes a stochastic process $\phi(t)=\phi(t,\omega)$  where $\omega$ is an element of a probability space, and similarly $n=n(t,\omega)$. 
The advantages of this probabilistic setting are: (i) in practical
problems one is typically interested in a range of possible missing data, rather than in some specific missing datum; (ii)
in data assimilation problems a probabilistic formulations is essential;       
and, most important, (iii) it is often possible to say something general and useful about the
probability density of missing data when one cannot say anything useful about 
any particular missing datum (an example will be given in the present paper).
A further discussion of the merits of this formulation will be given in the concluding section of this
paper.

In the present paper we assume in addition that enough information about the the statistics of $\phi$ is available (for specifics, see section 4). 
Of course one wishes to eventually solve problems where  
this information 
is not available, but looking at what happens when it is available 
is already instructive. 
We disregard here all sources of uncertainty other than underresolution,
except for a comment
in the conclusions section, though they too can be handled by the our formalism.

It is of course desirable to make $n(t)$ as small
as possible. In the present paper a low noise
intensity is reached by using the Mori-Zwanzig formalism of statistical mechanics \cite{Eva1,Fic,Zw4,Ch13,Darve}, which 
draws as much information as possible into the resolved part of the calculation.

In the next section we summarize the Mori-Zwanzig formalism in the version we need.
In section 3 we
present our example. In section 4 we present numerical results. Conclusions are drawn at the end of the paper.

There is a large literature on noise modeling. 
A general theoretical overview can be found in \cite{stuart}.
Related numerical work is summarized in \cite{Gha}.
Powerful techniques aimed at geophysical problems can be found in \cite{Maj2,Kal}.
A particularly interesting  approach, for situations where the problem can be
approximately viewed as linear, can be found in \cite{Miller}.
The work in the present paper has a kinship with the stochastic parametrization 
proposal of Wilks \cite{Wilks}, see also \cite{Hami}, 
though the mathematical and statistical tools differ,
and we are estimating the noise rather than parametrizing, 
i.e., estimating the variance rather than trying to improve 
the calculation of the mean. The technical difficulties discussed below
make us doubt that the latter goal is achievable, but if it is,
the Mori-Zwanzig formalism could be useful there as well. 

The use of the Mori-Zwanzig formalism for estimating uncertainty was pioneered in \cite{Stinis}
where a much more sophisticated implementation was used to quantify the uncertainty due to
uncertain parameters in a differential equation, without assuming that the full solution of the 
problem was already available. However, the simpler analysis in the present paper 
lends iself well to the particular purpose here, which to demonstrate the persistence of memory 
in underresolved dynamics. 

\section{The Mori-Zwanzig (MZ) formalism}

We largely follow here the exposition in \cite{Ch13} with suitable modifications
; see
also \cite{Ch8}. Note that our main tool, equation (\ref{newMZ}) below, can be derived 
without the full MZ formalism, as indicated there; we provide the more general 
formalism because it is needed for some comparisons later.

Consider again the system (\ref{main}) above,
\begin{equation}
\frac{d}{dt}\phi = R (\phi(t)),
\;\;
\phi(x,0)=x.
\;\;
\end{equation}
Denote the solution of this equation by $\phi(x,t)$, making its dependence on the
initial conditions explicit. 
Partition 
the vector $\phi$ as above 
into resolved variables $\hat{\phi}$ and unresolved variables $\tilde{\phi}$, and  
similarly set 
$x=(\hat{x},\tilde{x})$ and $R=(\hat{R},\tilde{R})$. 
Form the Liouville partial differential equation
\begin{equation}
u_t=Lu,
\label{liouville}
\end{equation}
where $u=u(x,t)$ satisfies the
initial condition $u(x,0)=g(x)$, $g$ is a given function and  
\begin{equation}
L=\sum R_j(x)\frac {\partial}{\partial x_j}. 
\label{liouville}
\end{equation}
One can show that the solution of this equation is $u(x,t)=g(\phi(x,t))$, where $\phi(x,t)$ is the solution of equations (\ref{main}). In particular,
if $g(x)=x_j$, it follows that $u(x,t)=\phi_j(x,t)$, where $x_j, \phi_j$ are the $j-$th components of respectively $x$ and $\phi$. 
Denote the solution of the Liouville equation with datum $g(x)$ by $u=e^{tL}g(x)$ (this is the ``semigroup notation"). In this notation, the identity $u(x,t)=g(\phi(x,t))$ becomes
\begin{equation}
e^{tL}g(x)=g(e^{tL}x).
\label{substitution}
\end{equation}
This identity shows that that the Liouville equation, which is linear, is equivalent to the (generally nonlinear) system of ordinary
differential equations (\ref{main}); the solution of equations (\ref{main}) provides the solution
of the Liouville equation, and conversely, if the Liouville equation can be solved for any initial
function $g(x)$, one obtains the $i-th$ component $\phi_i$ of (\ref{main}) by solving the Liouville 
equation with $g(x)=x_i$. 
One can also prove the identity:
\begin{equation}
e^{tL}L=Le^{tL}.
\label{commute}
\end{equation}

At time $t=0$ assign to the initial data $\hat{x}$ specific fixed values, and sample the initial data $\tilde{x}$ from a suitable probability density, as described in the introduction. 
The conditional expectation $E[h(x)|\hat{x}]$ of a function $h$ of $x$ given $\hat{x}$
is well-defined, and is an orthogonal projection of $h$ on the space of functions of $\hat{x}$.
Define the conditional expectation operator at the initial time $P$ by 
$Ph=E[h|\hat{x}]$ for any function $h$, and define and $Q=I-P$, where $I$ is the
 identity operator. Clearly $P,\,Q$ are both orthogonal projections and $P+Q=I$.

Pick as an initial datum for the Liouville equation the vector $\hat{x}$, i.e., consider only the
first $k$ equations in equation (\ref{main}). 
Using the notation $u=e^{tL}\hat{x}$, the Liouville equation 
becomes:
\begin{equation}
\frac {\partial}{\partial t}e^{tL}\hat{x}=Le^{tL}\hat{x}=e^{tL}L\hat{x}=e^{tL}(P+Q)L\hat{x},
\label{eq1}
\end{equation}
where the commutation relation (\ref{commute}) has been used. One can readily check
that $Lx_j=R_j$, so that $$PL\hat{x}=E[\hat{R}|\hat{x}],$$ which is a function of $\hat{x}$ only;
call this function $\bar{R}=\bar{R}(\hat{x})$. It follows that $e^{tL}QL\hat{x}=\hat{R}(\phi)-\bar{R}(\hat{\phi})$.
Equation (\ref{eq1}) is equivalent to the equation
\begin{equation}
\frac {d}{dt} \hat{\phi}=\bar{R}(\hat{\phi})+\left(\hat{R}(\phi)-\bar{R}(\hat{\phi})\right).
\label{newMZ}
\end{equation}
We propose to set $\mathcal{R}=\bar{R}$ in equation(\ref{psi0}) which defines 
the underresolved approximation, so that
\begin{equation}
\frac {d}{dt}\psi=\bar{R}(\psi).
\label{mm}
\end{equation}
The noise, as defined in equation (\ref{trick}), is
$n(t)=\hat{R}(\phi)-\bar{R}(\hat{\phi})$. 
Equation (\ref{newMZ}) is the main tool used in this paper; note
that one can get to it by simply setting $\mathcal{R}=\bar{R}$ in equation (\ref{psi}); the full MZ development is presented because it is
needed in the discussion.                

The derivation of the Mori-Zwanzig generalized Langevin equation  
requires several more steps.
Any two linear operators $A$ and $B$ satisfy the following identity
(the Dyson or Duhamel formula):
\begin{equation}
e^{t(A+B)}=e^{tA}+\int_0^te^{(t-s)(A+B)}Be^{sA} \, ds.
\end{equation}
Substituting 
$A=QL$ and $B=PL$, this becomes: 
\begin{equation}
\label{eq:dyson1}
e^{tL}=e^{tQL}+\int_0^t e^{(t-s)L}PLe^{sQL} \, ds.
\end{equation}
Substituting into (\ref{newMZ}), one finds: 
\begin{equation}
\frac{\partial}{\partial{t}} e^{tL}\hat{x}=
e^{tL}PL\hat{x}+e^{tQL}QL\hat{x}+
\int_0^t e^{(t-s)L}PLe^{sQL}QL\hat{x}\, ds,
\label{MZ}
\end{equation}
where $e^{tL}\hat{x}=\hat{\phi}$. 
This is the Mori-Zwanzig (MZ) generalized Langevin equation. 
The evaluation of the last two terms is in general difficult; 
an algorithm was proposed in \cite{Ch8} but it is too laborious for
practical use. Equation (\ref{MZ}) is in general not Markovian, and the
memory term is in general
non-zero, so that one cannot expect the noise $n(t)$ in equation (\ref{newMZ}) to have zero mean.

The MZ equation can be simplified by approximating 
the operator $e^{tQL}$ in the integrand of the third term
by $e^{tL}$ (see \cite{Ch8,Ch13}). 
With this approximation, the integral term simplifies to:
$$\int_0^t  e^{tL}PLQL\hat{x}ds = t e^{tL}LQL\hat{x}.$$ The memory term has been
reduced to a differential operator multiplied by the time $t$; the time starts at $t=0$ when the
initial values are assigned and when there is no uncertainty in the resolved variables. 

If one is interested only in the conditional expectations of $\hat{\phi}(t)$, one can
premultiply equation by the conditional averaging operator $P$; the noise term then drops out,
so that the MZ equation with the simplified integral term becomes
\begin{equation}
\frac {\partial}{\partial t}P\phi(t)=Pe^{tL}PL\hat{x}+tPe^{tL}PLQL\hat{x}.
\end{equation}
In the absence of the noise term the conditional expectation of a nonlinear function of $\phi$ cannot
be ascertained without approximation; we resort here to ``mean field" closure $PH(\phi)=H(P\phi)$
for any function $H=H(\phi)$, in particular for the functions $G(\phi)=PLQL(\hat{\phi})$ and $F(\phi)=PL(\hat{\phi})$, so that 
\begin{equation}
\frac {\partial}{\partial t}\Phi=F(\Phi)+tG(\Phi),
\label{tmodel}
\end{equation}
where $\Phi=P\phi$. 
This is the $t$-model approximation of the MZ equations (see \cite{Ch8,Ch13}).

\section{A model problem}

The dynamical system to which the theory of the preceding section will be applied is now presented. It is 
a Hamiltonian system with $m+1$ particles, one of which is resolved, or ``tagged", and the others are unresolved. The Hamiltonian
is
\begin{equation}
H=(1/2)(p^2+q^2+\sum_1^mp_i^2/g_i +(\sum_1^m q_i^2)/\epsilon+\alpha q^2\sum_1^mq_i^2),
\end{equation}
where $q,p$ are the position and the momentum of the tagged particle, $q_i,p_i$ are the positions 
and momenta of the $m$ other particles, $\alpha, \epsilon,$ and $g_i,i=1,\dots,m$ are parameters; $\epsilon$ quantifies the amount by which the unresolved particles are faster that the tagged particle; $\alpha$ quantifies the strength of the interaction between the tagged and unresolved particles, the parameters $g_i$ control the periods of the unresolved particles; in the examples below, $\alpha=\epsilon=m^{-1}$ and
$g_1=\sqrt{5}/2,\, g_2=\sqrt{3}/2.$ Except when stated otherwise, $m=2.$  
The 
underresolved system is one in which only the tagged particle is solved for, so that in the notations of the introduction, $k=2$ and $K=2m+2$. 
Note that  the notation for the variables has changed from $\phi$ to $q,p$, the latter being more transparent for a Hamiltonian system.  This system of equations is a generalization of
the example used in \cite{Ch8,Ch13} to demonstrate the properties of the MZ formalism; it resembles the model
problems studied in \cite{Ford,Zw3} but is more strongly nonlinear. This system is not ergodic; to demonstrate this,
figure 1 shows the projection of a single trajectory of the system on the $(q,p)$ plane; if the system were
ergodic, there would be no hole in the donut. 
\begin{figure}[h!tbp]
\begin{center}
{\includegraphics[width=0.75\textwidth]{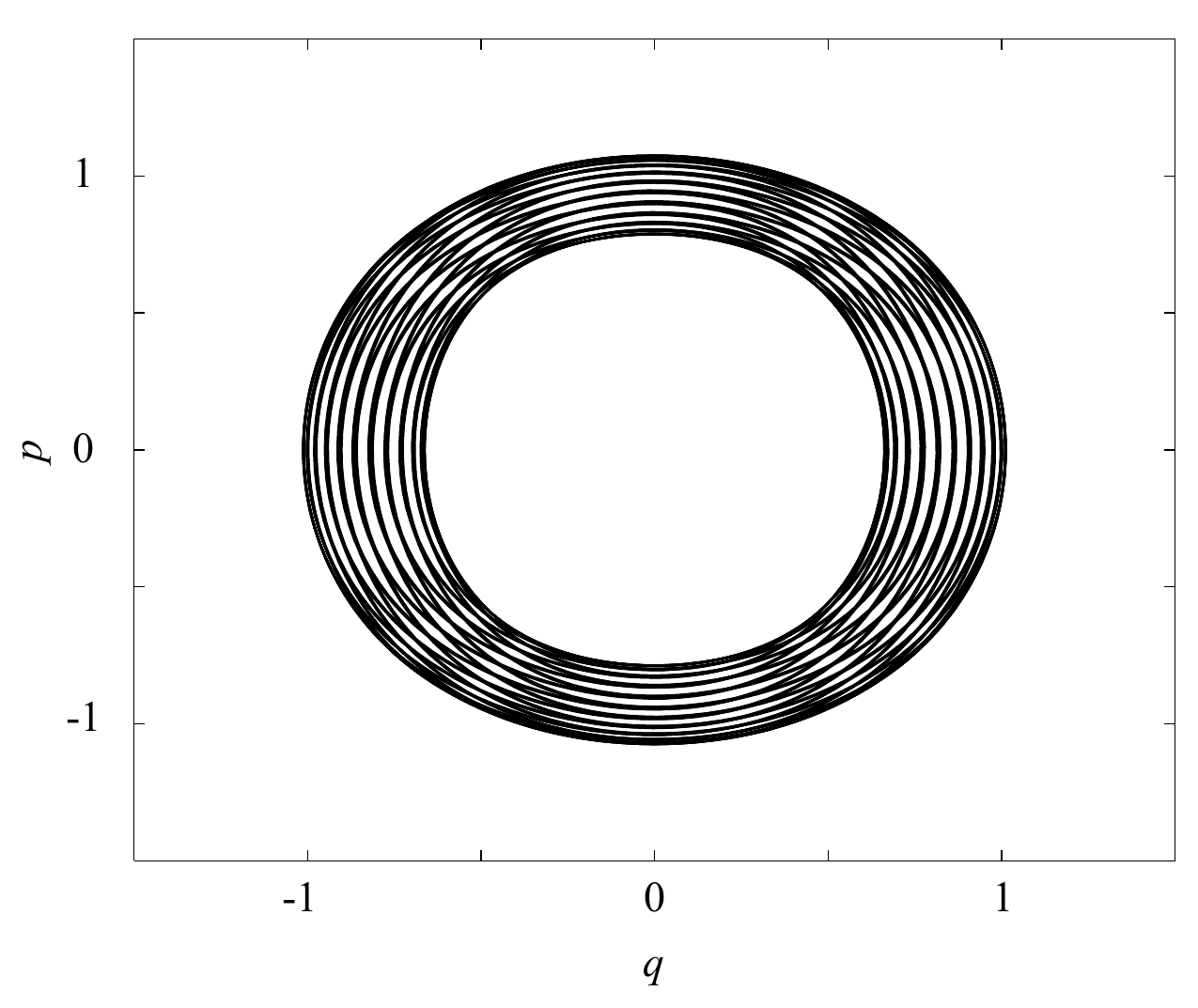}}
\caption{Ergodicity check.} 
\end{center}
\end{figure} 
A non-ergodic system is studied because of the expectation that
the constructions below will 
be used in problems with persistent states, as in weather forecasting (see e.g. \cite{Yang}).

A time $t=0$ one specifies initial values for $q(0),p(0)$ for the tagged variable, and one samples the $q_i,p_i$
from the canonical density $e^{-H/T}/Z$ conditioned by $q(0),p(0)$; here $T$ is a ``temperature" (which will be set equal to $1$ for simplicity) and $Z$ is the
partition function. The equations of motion are:
\begin{align}
\frac {d}{dt}q=&p, \nonumber\\
\frac {d}{dt}p=&-q(1+\alpha \sum_1^m q_i^2), \nonumber\\
\frac {d}{dt}q_i=&p_i/g_i, \,\,\,\,\,\, i=1,\dots,m, \nonumber\\
\frac {d}{dt}p_i=&-q_i(1+\epsilon \alpha q^2)/\epsilon, \,\,\,\,\, \,i=1,\dots,m.
\label{full}
\end{align}

The Liouville operator $L$ is:
\begin{equation}
L=p\frac {\partial}{\partial q}-q(1+\alpha \sum_1^m q_i^2)\frac {\partial}{\partial p}
+\sum (p_i/g_i)\frac {\partial}{\partial q_i}-\sum_1^m (q_i(1+\epsilon \alpha q^2)/\epsilon ) \frac {\partial}{\partial p_i}.
\end {equation}

Elementary calculations yield:
\begin{equation}
E[q(1+\alpha  \sum_1^m q_i^2|q,p]=q\left( 1+\frac {m\alpha \epsilon}{1+\alpha \epsilon q^2} \right),
\label{reducedp}
\end{equation}
while of course $E[p|q,p]=p$. The equations $\frac {d}{dt}\hat{\phi}=E[\hat{R}|\hat{\phi}]$
reduce to a two-component Hamiltonian system with Hamiltonian $$\hat{H}=(1/2)(p^2+q^2+m\log (1+\alpha \epsilon q^2)).$$
Similarly,
the vector $\hat{R}(\phi)-E[R(\phi)|\hat{\phi}]$ has components $(0,n(t))$ where 
\begin{equation}
n(t)=-q \left(\frac{m \alpha \epsilon}{1+\alpha \epsilon q^2}-\alpha \sum_1^mq_i^2 \right).
\label{forcing}
\end{equation}
In doing these calculations, it was assumed that:
$$E[f(q,q_i,p,p_i)|q,p]=\frac {\int f(q,q_i,p,p_i)e^{-H}\,dq_i\,dp_i}
{\int e^{-H}\,dq_i\,dp_i}$$
for an arbitrary function $f$, 
where $q_i$ stands for $q_1,q_2,\dots,q_m$, similarly for $p_i$; $dq_i$ stands
for $dq_1\cdots dq_m$, and similarly for $dp_i$; the Hamiltonian $H$ is 
a function of all the variables. The integrations are over the whole hyperplane
where $q$ and $p$ are constant; the non-ergodicity of the system is not taken
into account. 

Suppose one wants to calculate the conditional expectation of
the resolved variables $q,p$ given their initial values $q(0),p(0)$. In the low-dimensional example 
here 
this calculation can be done from first principles. 
Given $q,p$ at time $t=0$, one can repeatedly sample the coordinates and momenta of
the unresolved particles from the initial probability density  $e^{-H}/Z$; for each set of initial data, one can 
solve the full system of equations up to a time $t$; one can then average the values of 
$q(t),p(t)$ in 
these solutions; the result is, by definition, the conditional expectations $E[q(t)|q(0),p(0)], E[p(t)|q(0),p(0)]$.

On the other hand, one can derive and solve the MZ equations in the $t-$model
approximation. Elementary calculations show that 
these equations reduce to:
\begin{align}
\frac {d}{dt}Q=&P,\\
\frac {d}{dt}P=&-Q(1+\frac {m\alpha\epsilon}{1+\alpha \epsilon Q^2})
-\frac {2m \alpha^2 \epsilon^2 Q^2 P t}{(1+ \alpha \epsilon Q^2)^2},
\end{align}
where $Q=E[q(t)|q(0),p(0)], P=E[p(t)|q(0),p(0)].$ (For a related example
where the intermediate steps in the derivations are presented in detail, see \cite{Ch13}). 

In figure 2 the values of $Q=E[q(t)|q(0),p(0)]$ computed from
first principles, with no approximations other than use of the law of
large numbers, are compared with the results of $t$-model, with   
$m=2$.  
\begin{figure}[h!tbp]
\begin{center}
{\includegraphics[width=0.75\textwidth]{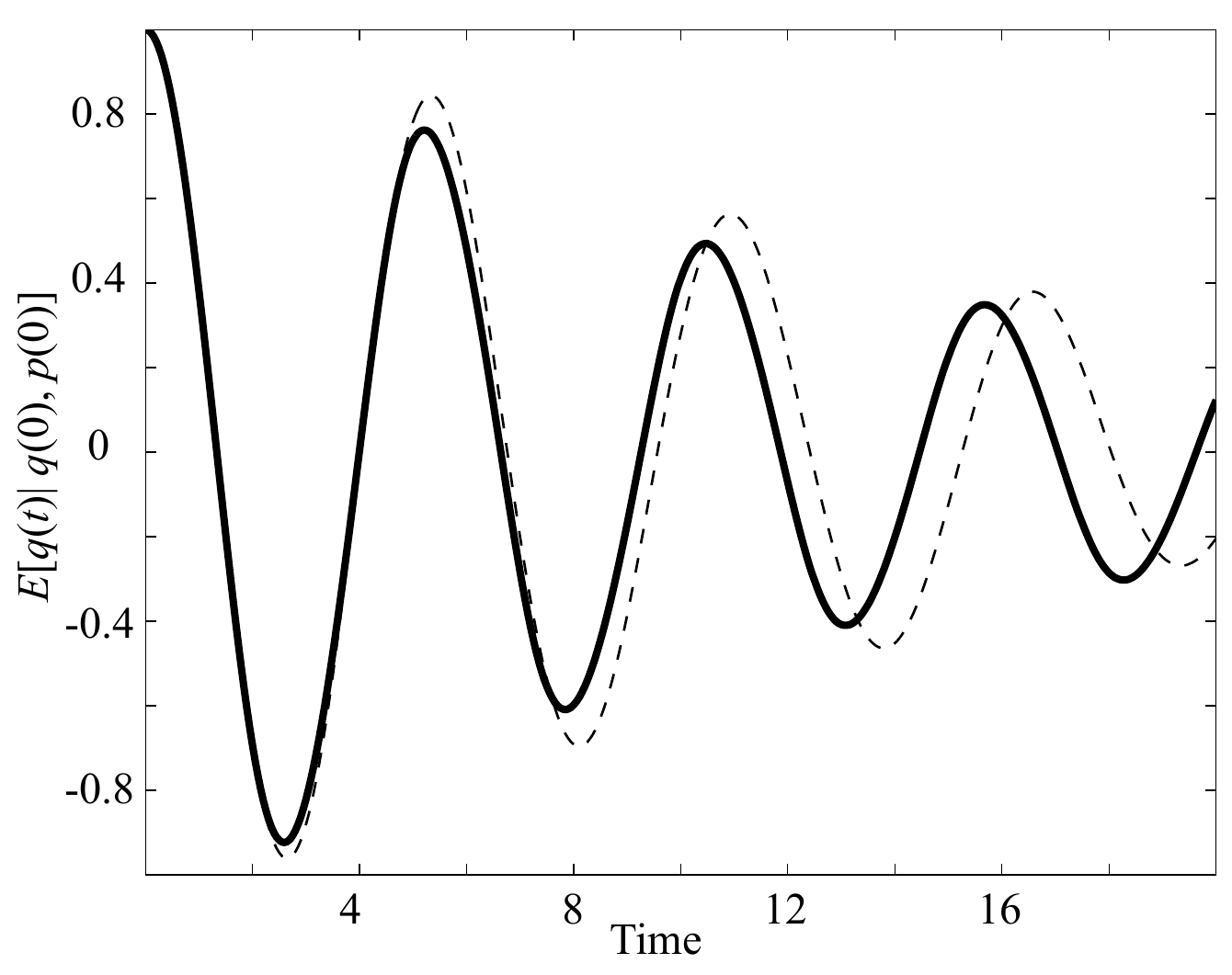}}
\caption{Time decay of conditional expectations: $t$-model (thick line) vs. the truth (thin dashed line).} 
\end{center}
\end{figure} 
It should be obvious
that the MZ formalism would give less noisy solutions when $m$, the number of unresolved
particles, increases, because then the instantaneous average of their
effects on the tagged particle approaches its average with respect
to the invariant measure, making the problem less interesting. We list again the approximations
used in deriving the $t-$model results: (i) the substitution $e^{tQL}=e^{tL}$ in the memory
term; (ii) the mean field approximation; (iii) no account taken of the
nonergodicity of the dynamics; in addition, (iv),
the measure with density $e^{-H}/Z$ is invariant at equilibrium, but since at $t=0$
the location and momentum $q,p$ of the tagged particle are fixed, the system
is not in equilibrium. One cannot expect the agreement in figure 2 to be
perfect. 

The conditional expectations decay in time because, in the presence of noise, the predictive
value of initial data decays in time. The uncertainty in the unresolved data affects the 
resolved variables, each random solution goes its own way, and while each preserves the Hamiltonian,
the average of the solutions converges to the equilibrium average, which here is zero. The decay 
in the expected value is due entirely to the uncertainty.

\section{Numerical analysis of the noise}

The program outlined in section 2 will now be implemented in the special case
described in the preceding section. The first step is the analysis of the
noise term $n(t)$ described in equation (\ref{forcing}). 
We have assumed that the noise is known in principle-- if the solutions of the full system
are known, formula (\ref{forcing}) above represents the noise. To be useful, for example
in data assimilation, one needs an explicit representation of the noise in terms
of a small number of parameters, and the goal in the present section is
to derive such a representation approximately.

The function $q(t)$ is resolved
and therefore available at every step, so the function one has to represent is: 
\begin{equation}
z(t)=n(t)/q(t)=- \frac{m \alpha \epsilon}{1+\alpha \epsilon q^2}+\alpha \sum_1^mq_i^2 ,
\label{noverq}
\end{equation} 

The system (\ref{main}) is stationary with an invariant canonical probability 
density 
$e^{-H}/Z$. Initial data (including data for the tagged variable) were sampled from this invariant density repeatedly; for each initial vector, the differential equations (\ref{full}) were solved numerically, the function $z(t)$ was evaluated, and the covariance function
$C(\tau)=E[z(0)z(\tau)]$ was computed. We observed that as $\tau \rightarrow \infty$, $C(\tau)$ converged to a constant $A$. We exhibit this convergence in
table 1, where we list values of $C(\tau)$ for various values of $\tau$, as
well as values of the integrated covariance $C^*(\tau)=\tau^{-1}\int_0^{\tau} C(s)ds$
which, in this oscillatory system, converge to their limit faster than $C(\tau)$ as $\tau$ increases.
\begin{table}[ht]
\caption{Convergence of $C(\tau)$ to a limit}
\centering
\begin{tabular}{c c c}
$\tau$&   $C(\tau)$  &       $C^*(\tau)$      \\ [0.5ex] 
\hline
     0 & 0.1967 & 0.1967\\
     1 & 0.0046 & 0.1018\\
     2 & 0.1771 & 0.0961\\
     3 & 0.0289 & 0.0999\\
     4 & 0.1398 & 0.0943\\
     5 & 0.0621 & 0.0971\\
     6 & 0.1072 & 0.0940\\
     7 & 0.0863 & 0.0949\\
     8 & 0.0909 & 0.0941\\
     9 & 0.0927 & 0.0935\\
    10 & 0.0922 & 0.0941\\
    11 & 0.0860 & 0.0928\\
    12 & 0.1009 & 0.0937\\
    13 & 0.0783 & 0.0926\\
    14 & 0.1054 & 0.0932\\
    15 & 0.0765 & 0.0926\\
    16 & 0.1062 & 0.0928\\
    17 & 0.0787 & 0.0926\\
    18 & 0.1041 & 0.0926\\
    19 & 0.0828 & 0.0927\\[1ex]
\hline
\end{tabular}
\end{table}
We then decomposed $C(\tau)$ in the form:
\begin{equation}
C(\tau)=\beta(\tau)+A,
\label{form}
\end{equation}
where $\beta(\tau)$ tends to zero as $\tau$ increases.
\begin{table}[ht]
\caption{Computed covariances}
\centering
\begin{tabular}{c c c}
$\tau$ & $C(\tau)$ & $\beta(\tau)$ \\ [0.5ex] 
\hline
  0.00 & 0.1967 & 0.1008 \\
  0.05 & 0.1953 & 0.0994 \\
  0.10 & 0.1917 & 0.0958 \\
  0.15 & 0.1860 & 0.0901 \\
  0.20 & 0.1782 & 0.0823 \\
  0.25 & 0.1686 & 0.0727 \\
  0.30 & 0.1574 & 0.0615 \\
  0.35 & 0.1449 & 0.0490 \\
  0.40 & 0.1313 & 0.0354 \\
  0.45 & 0.1170 &0.0212 \\
  0.50 & 0.1024 & 0.0065 \\[1ex]
\hline
\end{tabular}
\end{table}
Furthermore, for each set of initial data,   
the quantity $\tau^{-1}\int_0^{\tau}z(0)z(s)ds$, whose expected value is $C(\tau)$, individually converged to a limit, say $a$,
which is a function of the random initial values, so that $a=a(\omega)$,
with 
\begin{equation}
A=E[a],
\label{Aeqa2}
\end{equation}
where $A$ is the constant in equation (\ref{form}). The parameter $a$ is a ``persistence parameter",
embodying the fact that in a non-ergodic system the initial values are never fully forgotten.
It was pointed out in section 2 that one cannot expect $n(t)$ to have a zero mean.
The variable $a$ is roughly log-normal; in figure 3 we display a
histogram of $\log a$. 
\begin{figure}[h!tbp]
\begin{center}
{\includegraphics[width=0.8\textwidth]{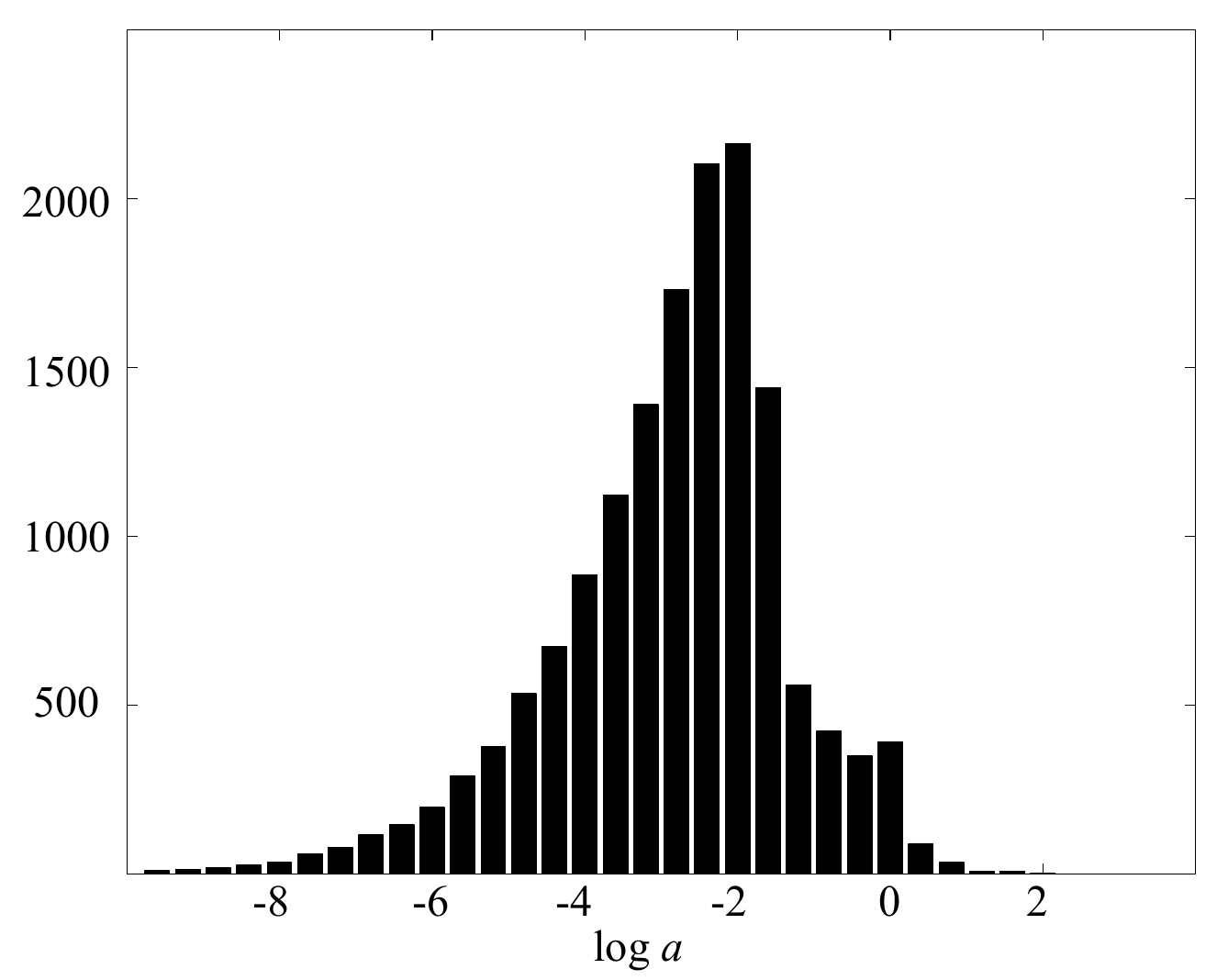}}
\caption{Histogram of $\log a$.} 
\end{center}
\end{figure} 
The calculations that produced these figures were performed
with a fourth-order Runge-Kutta method with time step $k=0.05$, repeated 20000 times to get the
statistics. 

The stochastic process $n(t)$ was modelled as the sum of a stationary time series $b(t,\omega)$, with a discrete step $k$ and with a covariance that approximates $\beta(\tau)$, to which, at each step, was added the square root of a  sample of the persistence parameter $a$,
sampled at time $t=0$ and left constant during each time evolution; this satisfies equation (\ref{form}). 
The covariance of the time series was approximated by a covariance
of the form $C(\ell k)=\beta_0\,d^{|\ell|}$, where $\ell$ is an integer, the coefficient $\beta_0$ is read from table 2, and the parameter $d$ is obtained by a least-squares fit of $C(\ell k)$ to the the computed function $\beta(\tau)$ in the interval $[0,\tau_0]$, where $\tau_0$ is the first zero of $\beta(\tau)$. 
In the case where $k=0.05$ and $n=2$, this gave the value $d=.902$; if one changes the numerical time
step from $k=0.05$ to $k=k^*$, $d$ becomes $d^*=d^{(k^*/0.05)}$. With this representation of $C(\ell k)$, the time series $b(\ell k)$ becomes a two-step recursion, so that if $b(jk)$ has been sampled
at time $jk$ and its value was $B(jk)$, then at time $(j+1)k$ the sample value is $B((j+1)k)=d\cdot B(jk)+\xi$,
where $\xi$ is sampled from an independent Gaussian distribution with mean zero and variance $\beta_0(1-d^2)$
(see e.g. \cite{Yagl,Ch13}). 
The values of $a$, one per  
sample of the initial data, were sampled using the histogram above.

The modeled stochastic process is the estimate of the uncertainty in the underresolved system we propose to use when an explicit representation of the noise is needed, for example in data assimilation.
Note that an equilibrium distribution of $n(t)$ is used though the system (\ref{full}) is solved with some of
the data having fixed values, i.e., it is not in equilibrium. This parallels the approximations made
in the application of the MZ formalism in section 2. 

As a check on the validity of this representation,
consider the solution of the system (\ref{psi}) with $n(t)$ sampled with the help of this representation. The system $\frac{d}{dt}\psi=\bar{R}(\psi)$ is Hamiltonian, so that the decay of the solutions of 
\begin{equation}
\frac{d}{dt}\psi=\bar{R}(\psi)+n(t), 
\label{bob}
\end{equation}
which 
should match the true decay, is due entirely to the presence of the noise $n(t)=p(t)z(t)$. In figure 4  we compare the average of the
solutions of equation (\ref{bob}), with fixed initial values but random forcing $n(t)$, to the exact
average computed in section 3.  
\begin{figure}[h!tbp]
\begin{center}
{\includegraphics[width=0.75\textwidth]{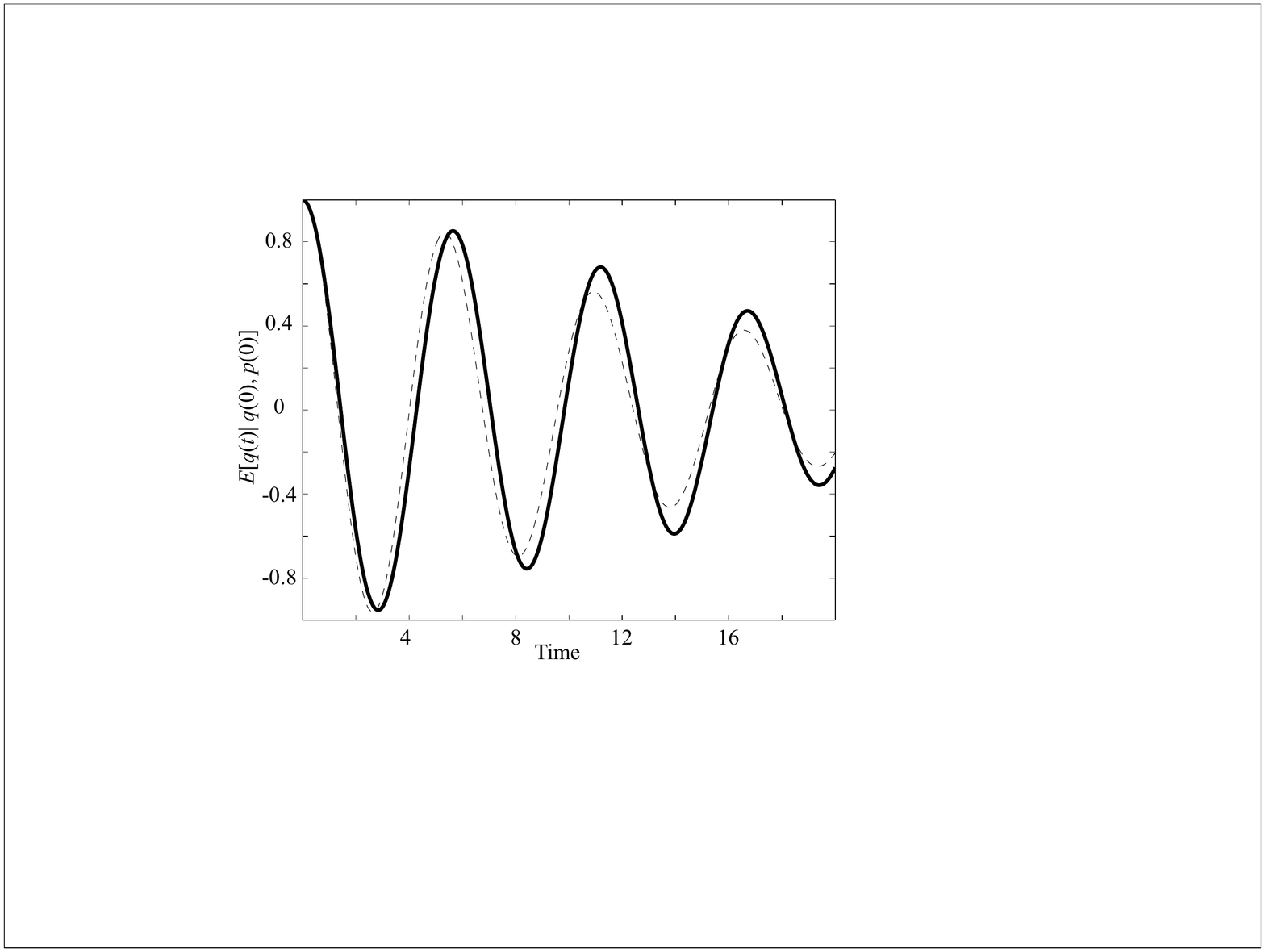}}
\caption{Time decay of conditional expectations: the noise model (thick line) vs. the truth (thin dashed line).} 
\end{center}
\end{figure} 

In figure 5 we compare the histogram of the values of $q(t)$ at time 
$t=20$ (a long time) obtained by solving repeatedly the full system (\ref{full}), with the
histogram obtained by solving the reduced system (\ref{bob}) where the
noise is sampled from the representation above. 
\begin{figure}[h!tbp]
\begin{center}
{\includegraphics[width=1\textwidth]{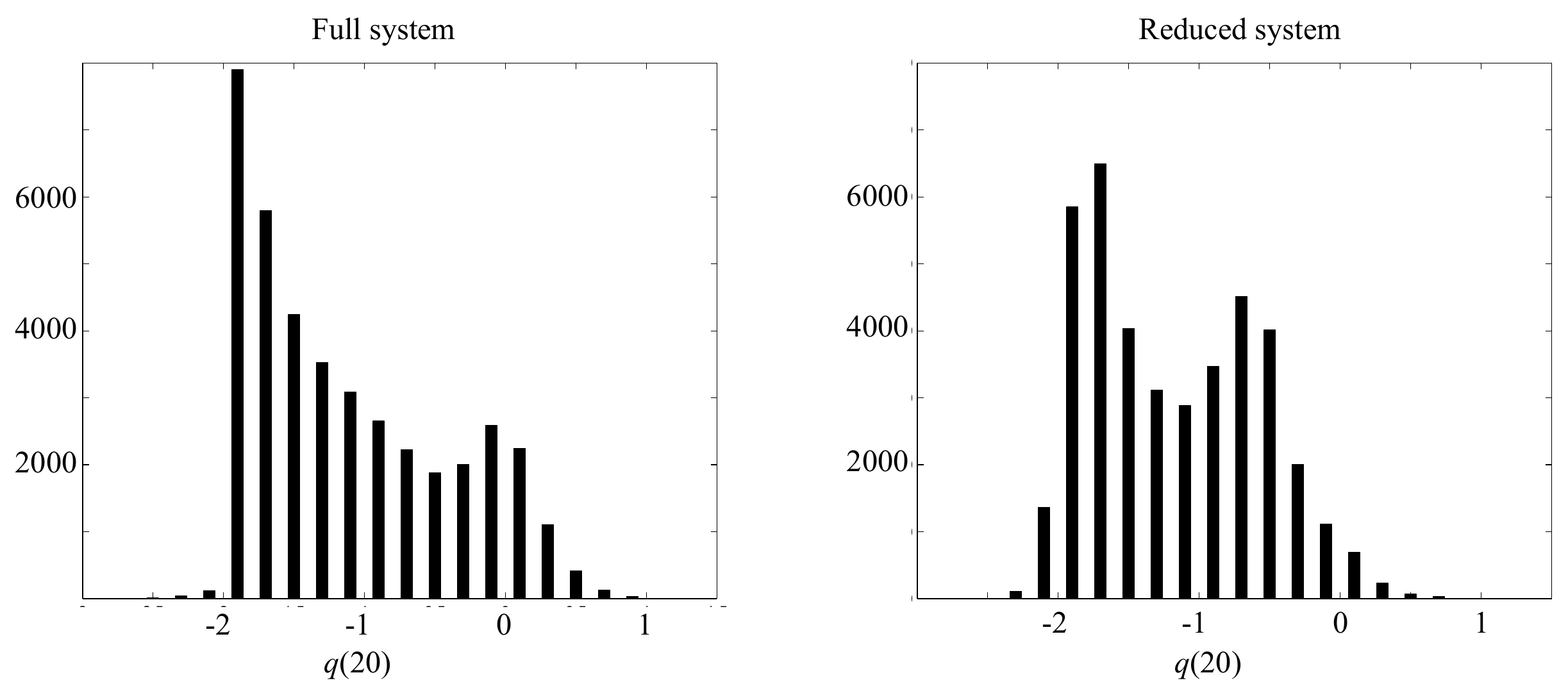}}
\caption{Histograms of q(t) at time t=20: the noise model vs. the truth.} 
\end{center}
\end{figure} 
The random equations (\ref{bob}) were solved by the
Klauder-Petersen scheme \cite{Kla1}; a rationale for this choice 
can be  found in \cite{Morz}. The time step was $k=0.01$. 
Note that equations (\ref{bob}) can be interpreted as a numerical implementation of
the full MZ Langevin equation (\ref{MZ}). 

These comparisons validate the representation of the
noise; indeed, the conditional expectations here are significantly closer to the truth than the one obtained with 
the $t$-model approximation.

\section{Conclusions}
In subsequent publications we expect to demonstrate the use of the formalism presented in this paper in more 
interesting cases. We think that the features of the noise seen in the particular case discussed above--
the fact that it is non-Gaussian and non-Markovian -- are general, and that 
the practice of data assimilation should take this into account. 

The non-Gaussianity is not surprising. The equations solved are non-linear and their invariant measure
is non-Gaussian, and there can be no expectation that the uncertainty in their solution is Gaussian.
The strongly non-Markovian nature of the noise is more interesting. It has been built into the model by the assumption
that the only randomness is in the the initial data, while the time evolution is deterministic. 
This is the right assumption for the assessment of the uncertainty due to underresolution in differential equations.
The fully resolved true solution of the problem satisfies a differential equation, one whose solution
one may be unable to find, but a differential equation none-the-less. The underresolved system one
can solve approximates a differential equation. The difference between solutions of differential
equations cannot be memoryless. These comments remain valid if the equations are chaotic- the size     
of the uncertainty may grow, but chaotic dynamics described by differential equations are not memoryless either.
The same is true if one chooses to model unresolved inputs as stochastic processes; 
they may modify the memory but not annihilate it. 

The present paper complements the work in \cite{CM13}, where it was shown that
in physically reasonable vector-valued data assimilation problems, the various one-time components of the noise must be
correlated; in particular, they must be spatially correlated if the several components describe 
physical variables estimated at different spatial points. In the present
 paper we show that the noise components must be are correlated in time as well.

\section{Acknowledgements}
We would like to thank Prof. R. Miller, of Oregon State University,
for advice, references, and
 illuminating discussions about error modelling in geophysical fluid dynamics, and our colleague Dr. Matthias Morzfeld for help with in the  
preparation of this paper. 
This work was supported in part by the Director,
Office of Science, Computational and Technology Research, U.S. Department
of Energy under Contract No. DE-AC02-05CH11231, and by the National Science Foundation
under grant DMS-1217065.


\begin{thebibliography}{10} 

\bibitem{Bens}
A. Bensoussan, Stochastic Navier-Stokes equations, Acta. Appl. Math. 38 (1995),
pp. 267-304.


\bibitem{Ch8}
A. J. Chorin, O.H. Hald, and R. Kupferman,
Optimal prediction with memory.
Physica D 166 (2002), pp. 239-257. 


\bibitem{Ch13}
A.J. Chorin and O.H. Hald, 
Stochastic Tools for Mathematics and Science. Springer-
Verlag, New York, 2013.  

\bibitem{Ch3}
A.J. Chorin and X. Tu, Implicit sampling for particle filters, Proc. Nat. Acad.
Sc. USA
106 (2009), pp. 17249-17254.


\bibitem{CM13}
A.J. Chorin and M. Morzfeld, Conditions for successful data assimilation,
 under review, J. Geophys. Res. 2013.

\bibitem{Darve}
A.J. Darve, J. Salomon, and A. Kia, Computing generalized Langevin equations and
generalized Fokker-Planck equations, Proc. Nat. Acad. Sc. USA 27 (2009), 
pp. 10884-10889.

\bibitem{Douc}
A. Doucet, N. de Freitas, and N. Gordon, Sequential Monte Carlo Methods in Practice,
Springer, NY, 2001.

\bibitem{Eva1}
D. Evans and G. Morriss, Statistical Mechanics of Nonequilibrium Liquids,
Academic Press, London, 1990.


\bibitem{Fic}
E. Fick and G. Sauerman, The Quantum Statistics of Dynamical Processes,
Springer, 1990. 

\bibitem{Ford}
G. Ford, M. Kac, and P. Mazur, Statistical mechanics
of assemblies of coupled oscillators, J. Math. Phys. 6 (1965), pp. 504-515.

\bibitem{Gha}
R. Ghanem and P. Spanos, Stochastic Finite Elements: A Spectral Approach,
Springer, New York, 1991.

\bibitem{Hami}
T. Hamill and J. Whitaker, Accounting for the error due to unresolved
scales in ensemble data assimilation: a comparison of different approaches,
Month. Weather Rev. 133 (2005), pp. 3132-3247.

\bibitem{Kal}
E. Kalnay, Atmospheric Modeling, Data Assimilation and Predictability,
Cambridge, 2003. 

\bibitem{Kla1}
J. Klauder and W. Petersen, Numerical integration of multiplicative-noise
stochastic differential equations, SIAM J. Num. Anal. 22 (1985), pp. 1153-1166.



\bibitem{Kuksin}
S. Kuksin and A. Shirikyan, Mathematics of Two-Dimensional Turbulence,
Cambridge, 2012.

\bibitem{Maj2}
A. Majda and J. Harlim, Filtering Complex Turbulent Systems, Cambridge, 2012. 

\bibitem{Miller}
J. Richman, R. Miller, Y. Spitz, Error estimates for assimilation of satellite
sea surface temperature data in ocean climate models, Geophys. Res. Lett.
32 (2005), DOI:10.1029/2005GL023591.


\bibitem{Morz}
M. Morzfeld, X. Tu, E. Atkins, and A.J. Chorin, A random map implementation of implicit 
filters, J. Comp. Phys. 231 (2012),pp. 2049-2066. 
    

\bibitem{Papa}
G. Papanicolaou, 
Introduction to the asymptotic
analysis of stochastic equations, in ``Modern Modeling of Continuum Phenomena",
R. DiPrima (ed.), Providence RI, 1974.

\bibitem{Stinis} Mori-Zwanzig reduced models for uncertainty quantification I: 
parametric uncertainty, arXiv:1211.4285v1 [math.NA] 19Nov2012.

\bibitem{stuart}
A.M. Stuart. Inverse problems: a Bayesian perspective,
Acta Numerica (19), 2010. pp 451-559

\bibitem{Wilks}
D. Wilks, Effects of stochastic parametrization in the Lorenz '96 system,
Quart. J. Roy. Meteo. Soc. 131 (2005), pp. 389-407.

\bibitem{Yagl}
A. Yaglom, An Introduction to the Theory of Stationary Random Functions,
Dover, New York, 1962.

\bibitem{Yang}
X. Yang, W. Ren, and E. Vanden Eijnden, Nonequilibrium statistical mechanics
of climate variability, submitted for publication, 2013.


\bibitem{Zw3}     
R. Zwanzig, Nonlinear generalized Langevin equations, J. Stat. Phys., 9, (1973),
pp. 215-220.

\bibitem{Zw4}
R. Zwanzig, Nonequilibrium Statistical Mechanics, Oxford, 2001.

\end{thebibliography}
\end{document}